\documentclass[a4paper, 10pt]{amsart}

\usepackage{geometry}                
\usepackage{graphicx}
\usepackage{amssymb}
\usepackage{epstopdf}
\usepackage{hyperref}
\usepackage{verbatim}
\usepackage[arrow, matrix, curve]{xy}

\newtheoremstyle{thms}
	{}{}{\itshape}{}{\bfseries }{}{ }
	{\thmname{#1} \thmnumber{#2}. \thmnote{\bfseries{[#3]}}}
\newtheoremstyle{name}
	{}{}{\itshape}{}{\bfseries }{}{ }
	{\thmname{#1}\thmnumber{#2}\thmnote{\bfseries{[#3]}}}
\newtheoremstyle{defs}
	{}{12pt}{\normalfont}{}{\bfseries }{}{ }
	{\thmname{#1} \thmnumber{#2}. \thmnote{\bfseries{(#3)}}}
\newtheoremstyle{rmk}
	{}{}{\normalfont}{}{\itshape }{}{ }{\thmname{#1}. \thmnote{#3}}
\newtheoremstyle{claim}
	{}{}{\normalfont}{}{\itshape}{}{ }{\thmname{#1} \thmnumber{#2}. \thmnote{#3}}

\theoremstyle{thms}
\newtheorem{Prop}{Proposition}[section]
\newtheorem*{Thm}{Theorem}
\newtheorem{Lemma}[Prop]{Lemma}

\theoremstyle{name}

\theoremstyle{defs}
\newtheorem{Def}[Prop]{Definition}

\theoremstyle{rmk}
\newtheorem*{Rmk}{Remark}

\theoremstyle{rmk}
\newtheorem*{warn}{Warning}

\theoremstyle{rmk}
\newtheorem*{ex}{Example}

\theoremstyle{claim}

\newenvironment{pf}{\begin{proof}[Proof]}{\end{proof}}

\newcommand{\R}{\mathbb{R}}

\newcommand{\C}{\mathbb{C}}

\newcommand{\h}{\mathcal{H}}
\newcommand{\ad}{\mathcal{L}}
\newcommand{\adc}{\mathcal{K}}

\newcommand{\F}{\mathcal{F}}

\newcommand{\e}{\mathcal{E}}

\newcommand{\adl}{\mathrm{ad}_\lambda}

\newcommand{\id}{\mathrm{id}}

\newcommand{\E}{\underline{E}}

\newcommand{\BL}{\mathcal{L}_{BG}}

\newcommand{\ind}{\mathrm{ind}}

\newcommand{\MF}{\mathrm{MF}}
\newcommand{\Or}{\mathrm{Or}}
\newcommand{\GJ}{\mathrm{GJ}}

\DeclareMathOperator*{\colim}{colim}

\title{The Analytical Assembly Map and Index Theory}
\author{Markus Land}
\address{Rheinische Friedrich-Wilhelms Universit\"at Bonn, Mathematisches Institut, Endenicher Allee 60, 53115 Bonn, Germany}
\email{land@math.uni-bonn.de}
\date{February 2014}

\begin{document}

\begin{abstract}
In this paper we study the index theoretic interpretation of the analytical assembly map that appears in the Baum-Connes conjecture. In its general form it may be constructed using Kasparov's equivariant $KK$-theory. In the special case of a torsionfree group the domain simplifies to the usual $K$-homology of the classifying space $BG$ of $G$ and it is frequently used that in this case the analytical assembly map is given by assigning to an operator an equivariant index. We give a precise formulation of this statement and prove it.

\end{abstract}

\maketitle

\tableofcontents

\begin{section}{Introduction}
\noindent
Let $G$ be a countable discrete group. The Baum-Connes conjecture predicts a certain analytical 
assembly map
\[\xymatrix{RK^G_*(\underline{E}G) \ar[r] & K_*(C^*_rG) }\]
to be an isomorphism. In the case where the group is torsionfree the domain of this map may be identified with the compactly supported analytic $K$-homology of the classifying space $BG$ of $G$ and it is frequently used that the resulting map
\[\xymatrix{RK_*(BG) \ar[r] & K_*(C^*_rG)}\]
is given by associating to an elliptic differential operator an equivariant index.
The goal of this paper is to prove that this is indeed the case.

\noindent
This is important for the following reason. One standard method of proving the Baum-Connes conjecture is by the socalled Dirac-dual-Dirac method. This uses the construction of the assembly map as proposed by Kasparov. But when one wants to prove that the Baum-Connes conjecture implies (for example) the trace conjecture (which in turn implies the Kaplansky conjecture) one uses the interpretation of the assembly map as a Mishchenko Index.
So when relating the Baum-Connes conjecture to other classical conjectures one needs the index theoretic interpretation of the analytical assembly map.

\noindent
To the author's knowledge, this result has not been published yet and is in fact more subtle than we expected. The main obstacle is to relate the Mishchenko bundle to the canonical projection associated to any proper and cocompact $G$-space, which is used in the Kasparov picture.

\noindent
Most arguments of the paper are containend in detail in the author's master's thesis, but a crucial step was still missing there. We close this gap using a recent result of Buss-Echterhoff about fixed point Hilbert-modules in a 
special case.

\noindent
The paper is devided into three parts.

\noindent
In section 2 we recall the definition of the (full) analytical assembly map as proposed by Kasparov.

\noindent
In section 3, we define an index map
\[\xymatrix{KK_*(C_0(BG),\C) \ar[r]^-{\MF} & KK_*(\C,C^*G)}\]
and relate it to classical Mishchenko-Fomenko index theory.

\noindent
In section 4 we give a proof that these constructions coincide. This implies that the same is true for the reduced versions. On the way we prove a factorization of the descent homomorphism as it appears in the analytical assembly map and give some recollections on the Morita theory needed for the proof.
\addtocontents{toc}{\protect\setcounter{tocdepth}{1}}
\subsection*{Acknowledgements}
I want to thank Wolfgang L\"uck for introducing me to the field of isomorphism conjectures, and especially the Baum-Connes conjecture. Without his support during the writing of my master's thesis, this paper would not have been possible. I want to thank Alain Valette for encouraging me to think about the index theoretic interpretation of the analytical assembly map and Nigel Higson for telling me about the factorization of the descent homomorphism as I need it. I am particularly indebted to Siegfried Echterhoff who answered a number of questions concerning $C^*$-algebraic methods used in the proof of my comparison theorem. I also want to thank the referee for helpful suggestions. This work has been supported by the Leibniz Preis of Wolfgang L\"uck.
\addtocontents{tox}{\protect\setcounter{tocdepth}{2}}
\end{section}

\begin{section}{The Kasparov Approach to Analytical Assembly}
\noindent
The main input in Kasparov's definition of the analytical assembly map are a descent homomorphism and a canonical element in $KK(\C,C_0(X)\rtimes G)$ for any proper and cocompact $G$-space $X$. We will first recall the descent homomorphism to fix notation.

\begin{Lemma}\label{descent homomorphism}\label{descent}
For any $G$-$C^*$-algebras $A$ and $B$ there is a descent homomorphism
\[\xymatrix@C=2.5cm{ KK_*^G(A,B) \ar[r]^-{j^G_{(r)}} & KK_*(A\rtimes_{(r)} G, B\rtimes_{(r)} G) }\]
which is functorial and compatible with Kasparov products in the obvious sense.
\end{Lemma}
\begin{pf}
This is due Kasparov \cite[Theorem 3.11]{KA} and is also explained in \cite[2.2]{BE}.
To fix notation let us briefly summarize \cite[2.2]{BE}.

We consider an equivariant $KK$-cycle given by $\lbrack \e,\pi,\F \rbrack \in KK^G_*(A,B)$ i.e., $\e$ is a $G$-Hilbert-$B$-module, $\pi:A \to \ad(\e)$ is a graded equivariant $^*$-homomorphism and $\F \in \ad(\e)$ an odd self-adjoint operator satisfying the usual compatibility relations.

We then consider $C_c(G,\e)$ as a pre-Hilbert-$C_c(G,B)$-module as in \cite[2.2]{BE}. There is a left action of $C_c(G,A)$ on $C_c(G,\e)$ using the $G$-action on $\e$. Now we can complete this to 
\[ \e\rtimes G = \overline{C_c(G,\e)} \]
which is then a Hilbert-$B\rtimes G$-module. The action of $C_c(G,A)$ on $C_c(G,\e)$ extends to a graded $^*$-homomorphism $\xymatrix{\tilde{\pi}:A\rtimes G \ar[r] & \ad(\e\rtimes G)}$.
Furthermore we define $\tilde{\F} \in \ad(\e\rtimes G)$ by $\tilde{\F}(\alpha)(g) = \F(\alpha(g))$ for $\alpha \in C_c(G,\e)$. Then
\[ j^G\lbrack \e,\pi,\F \rbrack = \lbrack \e\rtimes G,\tilde{\pi},\tilde{\F} \rbrack.\]
For the reduced descent homomorphism, we simply complete $C_c(G,\e)$ to a Hilbert-$B\rtimes_r G$-module $\e\rtimes_r G$ and it follows (using e.g. \cite[Lemma 3.9]{KA}) that the canonical morphism $A \rtimes G \to \ad(\e\rtimes_r G)$ factors over $A\rtimes_r G$ as needed.
\end{pf}

Now let $X$ be a proper and cocompact $G$-space. 

\begin{Lemma}\label{key for projection}
There exists a non-negative function $\psi \in C_c(X)$ such that $\sum_{g\in G} \psi(gx) = 1$ for all $x \in X$. Such a $\psi$ will be referred to as \emph{cut-off function}. Using this we define an element $p_X \in C_c(G\times X)$ by
\[ p_X(g,x) = \sqrt{\psi(x)\cdot\psi(g^{-1}x)}\]
and view it as an element $p_X \in C_0(X)\rtimes G$.
This element is a projection and its $KK$-theory class $\lbrack \,p_X \rbrack \in KK(\C,C_0(X)\rtimes G)$ is independent of the choice of the cut-off function.
\end{Lemma}
\begin{pf}
This is proven in the more general setting for groupoids in \cite[sections 6.2 and 6.3]{JLT}. Concrete calculations may also be found in \cite[section 2.3]{ML}.
\end{pf}

\begin{Def}
Let $X$ be a proper and cocompact $G$-space. Then we define a morphism $\mu_X$ by the composite
\[\begin{xy}
\xymatrix@C=1.3cm{ KK^G_*(C_0(X),\C) \ar[r]^-{j^G} & KK_*(C_0(X)\rtimes G, C^* G) \ar[r]^-{-\circ \lbrack\, p_X \rbrack} & KK_*(\C,C^*G) \cong K_*(C^* G).}
\end{xy}\]
\end{Def}

\begin{Rmk}
This map is natural with respect to equivariant maps of proper $G$-spaces because the canonical $KK$-class $\lbrack\, p_X \rbrack \in KK(\C,C_0(X)\rtimes G)$ does not depend on a specific $\psi$ as in Lemma \ref{key for projection}.
\end{Rmk}

\begin{Def}\label{analytic assembly}
For a countable discrete group $G$ the full \emph{analytical assembly map} is the map
\[\xymatrix{RK^G_*(\underline{E}G) \stackrel{\mathrm{Def}}{=} \colim\limits_{X\subset \underline{E}G} KK^G_*(C_0(X),\C) \ar[r]^-{ \mathcal{A}} & K_*(C^*G),}\]
where the colimit runs over all cocompact $G$-invariant subsets $X$ of $\underline{E}G$, the classifying space for proper $G$-actions, and the map is induced by the maps $\mu_X$.
By the previous remark this is well-defined. If $G$ is torsionfree then $\E G = EG$.
\end{Def}

\begin{Rmk}
The reduced analytical assembly map may be obtained from the previous full version by post composing with the canonical morphism $K_*(C^*G) \to K_*(C^*_r G)$.
\end{Rmk}

\end{section}

\begin{section}{The Mishchenko-Fomenko Index}

\begin{Def}
For a $CW$-complex $X$ and a unital $C^*$-algebra $A$ denote by $K(X;A)$ the Grothendieck group of the monoid of isomorphism classes of finitely generated projective Hilbert-$A$-module bundles over $X$ under direct sum.
\end{Def}

\begin{Prop}\label{A-bundles and K theory}
If $X$ is compact, there is an isomorphism 
\[\xymatrix{K(X;A) \ar[r] & KK(\C,C(X)\otimes A)}\] 
induced by assigning to such a finitely generated projective Hilbert-$A$-module bundle its module of sections.
\end{Prop}
\begin{pf}
This is proven in \cite[Proposition 3.17]{TS}.
\end{pf}

\begin{Def}
For a compact space $X$ and a map $\xymatrix{f: X \ar[r] & BG}$ we consider the $G$-bundle $\hat{X} \to X$ classified by $f$ and define 
the \emph{Mishchenko line bundle} to be the following associated bundle
\[ \mathcal{L}_{f} = \hat{X}\times_G C^*G,\]
where the action of $G$ on $C^*G$ is given by left multiplication.
If $f = \id: BG \to BG$ is the identity or the inclusion of a subspace $Y$ we simply write $\BL$ or $\mathcal{L}_Y$ for this bundle.
\end{Def}

\begin{Rmk}
Of course, by construction we have that
\[\mathcal{L}_f = f^*(\BL)\]
and by Proposition \ref{A-bundles and K theory} we view this element as
\[ \lbrack \mathcal{L}_f \rbrack \in KK(\C,C(X)\otimes C^*G).\]
\end{Rmk}

\begin{Def}
Now let $X \subset BG$ be a compact subset. We define a Mishchenko-Fomenko index map $\MF$ by the composite
\[\begin{xy} 
\xymatrix@C=1.8cm{KK_*(C(X),\C) \ar[r]^-{\tau_{C^*G}} & KK_*(C(X)\otimes C^* G, C^* G) \ar[r]^-{-\circ \lbrack \mathcal{L}_X \rbrack} & KK_*(\C,C^* G) }
\end{xy}\]
where $\mathcal{L}_X$ is the bundle $\BL$ restricted to the subset $X$. This is just the cup-cap product map with the element $\lbrack \mathcal{L}_X \rbrack \in KK(\C,C(X)\otimes C^*G)$
This construction induces a map on colimits :
\[\xymatrix@!C{RKK_*(C_0(BG),\C) \ar[r]^-{\MF} & KK_*(\C,C^* G)}\]
because $\mathcal{L}_X$ is natural in $X$ as it is the pullback of a bundle over $BG$.
\end{Def}

\begin{Rmk}
In \cite{BHS} it is shown that every element in $RKK(C_0(BG),\C)$ may be represented by a triple $\lbrack M,f,E\rbrack$ where $M$ is a $\mathrm{spin}^c$-manifold, $f:M \to BG$ is a continuous map and $E$ is a hermitian vector bundle over $M$. Using this we can relate the previous map $\MF$ to the classical construction of $C^*$-algebra valued indices as in \cite{MF}.
\end{Rmk}

\begin{Prop}
The index map
\[\xymatrix@!C@R=.2cm{RKK_0(C_0(BG),\C) \ar[r]^-{\MF} & KK_0(\C,C^*G) \\
			\lbrack M,f,E \rbrack \ar@{|->}[r] & \ind(D_E\rtimes \mathcal{L}_f)}\]
maps the class $\lbrack M,f,E \rbrack$ to the Mishchenko-Fomenko Index of the Dirac operator twisted by the Mishchenko line bundle as constructed in \cite{MF}.
\end{Prop}
\begin{pf}
This follows from \cite[Theorem 6.22]{TS} and the commutativity of the diagram
\[\begin{xy}
\xymatrix{ RK_0(BG) \ar[r]^-{\tau_{C^* G}} & RKK_0(C_0(BG)\otimes C^* G, C^* G) \ar[r]^-{-\circ \lbrack \mathcal{L}_{BG} \rbrack} & KK_0(\C,C^* G) \\
	 K_0(M)\ar[r]_-{\tau_{C^* G}} \ar[u]^{f_*} & KK_0(C(M)\otimes C^* G,C^* G) \ar[r]_-{-\circ\lbrack \mathcal{L}_f \rbrack} \ar[u]_{(f^*\otimes C^* G)^*} & KK_0(\C,C^* G). \ar@{=}[u]}
\end{xy}\]
The analogous statement for the reduced version (using $C^*_r G$ instead of $C^*G$) is also true by essentially the same reasoning.
\end{pf}

\end{section}

\begin{section}{The Comparison Theorem}
In this section we study the relationship between the analytical assembly map of section 2 and the Mishchenko-Fomenko index map of section 3 in the case of a torsionfree group.
We recall that these are maps
\[\xymatrix@!C{\MF : RKK_*(C_0(BG),\C) \ar[r] & KK_*(\C,C^*G)}\text{ and }\]
\[\xymatrix@!C{\mathcal{A} : RKK^G_*(C_0(E G),\C) \ar[r] & KK_*(\C,C^*G).} \]
The main result of this paper is the following
\begin{Thm}\label{main thm}
Let $G$ be a countable, torsion-free, discrete group. Then there is an identification of the domains such that the following diagram commutes:
\[\xymatrix@!C{RKK_*(C_0(BG),\C) \ar[r]^-{\MF} & KK_*(\C,C^*G) \\ RKK^G_*(C_0(E G),\C) \ar@/_1.2pc/[ur]_-{\mathcal{A}} \ar[u]^\cong&}\]
\end{Thm}

\noindent
So let us begin by explaining the identification
\[\xymatrix{RKK^G_*(C_0(EG),\C) \ar[r]^-{\cong} & RKK_*(C_0(BG),\C).}\]
This proceeds in two steps.

First, we use the dual of the Green--Julg Theorem as stated e.g. in \cite[20.2.7 (b)]{BB}.
\begin{Prop}
For a discrete group $G$ and a $G$-$C^*$-algebra $A$ there is a canonical isomorphism
\[\xymatrix@!C{\GJ : KK_*^G(A,\C) \ar[r]^-{\cong} & KK_*(A\rtimes G,\C)}\]
i.e. the equivariant analytical $K$-homology coincides with the unequivariant $K$-homology of the full crossed product.
\end{Prop}

The next result we need is the following proposition due to Green, see \cite{PGr}.
\begin{Prop}
Suppose that a discrete group $G$ acts properly and freely on a space $X$, e.g. the action is proper and the group is torsionfree. Then the algebras $C_0(X)\rtimes G$ and $C_0(X/G)$ are Morita equivalent.
\end{Prop}
\begin{pf}
There are at least two ways to construct this Morita equivalence, and since we will need both descriptions we  briefly mention both.
Both construct an imprimitivity $C_0(X/G)$-$C_0(X)\rtimes G$ bimodule.
A canonical way is to consider the module
\[ \F_c(X) = C_c(X) \]
with bimodule structure given by
\begin{enumerate}
\item $(f.\varphi)(x) = f\lbrack x \rbrack \cdot \varphi(x)$, for $f \in C_0(X/G)$ and $\varphi \in C_c(X)$,
\item $(\varphi.\alpha)(x) = \sum\limits_{g \in G} \varphi(g^{-1}x)\cdot \alpha(g^{-1},g^{-1}x)$, for $\alpha \in C_c(G\times X)$,
\end{enumerate}
and inner products given by
\begin{enumerate}
\item $\langle \varphi,\varphi^\prime \rangle_{C_0(X)\rtimes G}(g,x) = \overline{\varphi(x)}\cdot \varphi^\prime(g^{-1}x)$, as well as
\item $_{C_0(X/G)}\langle \varphi,\varphi^\prime \rangle \lbrack x \rbrack = \sum\limits_{g\in G} \overline{\varphi(g^{-1}x)}\cdot \varphi^\prime(g^{-1}x)$.
\end{enumerate}
This completes to a Hilbert-$C_0(X)\rtimes G$-module which we call $\F(X)$. It still carries the structure of a $C_0(X/G)$-$C_0(X)\rtimes G$-bimodule and $C_0(X/G) \cong \adc(\F(X))$. Moreover the $C_0(X)\rtimes G$-valued inner product is full if and only if the $G$-action on $X$ is free, and so this is an imprimitivity bimodule as needed.

The other approach uses the projection $p_X \in C_0(X)\rtimes G$ and general Morita theory associated to projections.
It is a general fact about corners that given any $C^*$-algebra $A$ and a projection $p \in A$ the module $pA$ with the obvious structure becomes an imprimitivity $pAp$-$\overline{ApA}$ bimodule. A projection is called \emph{full} if $\overline{ApA} = A$ and if $G$ acts freely on $X$ then the projection $p_X$ is full in this sense.
Hence $p_X\cdot (C_0(X)\rtimes G)$ is also an imprimitivity bimodule as stated in the proposition since the corner $p_X\cdot (C_0(X)\rtimes G)\cdot p_X$ is isomorphic to $C_0(X/G)$, see \cite{EE}.
But the projection is not canonical (only its $K$-theory class is canonical) so this is a draw-back in this definition. We have an (noncanonical) isomorphism of imprimitivity bimodules
\[\xymatrix@!C{\Phi: \F(X) \ar[r] & p_X\cdot(C_0(X)\rtimes G)}\]
which restricted to $\F_c(X)$ is given by
\[ \Phi(\varphi) = \langle \Theta, \varphi\rangle_{C_0(X)\rtimes G} \]
where $\Theta = \sqrt{\psi}$ for a cut-off function $\psi$ as in Lemma \ref{key for projection}.
\end{pf}

\begin{Def}
We denote the resulting (invertible) $KK$-element by
\[ \lbrack \F(X) \rbrack \in KK(C_0(X/G), C_0(X)\rtimes G).\]
\end{Def}
\begin{Rmk}
The element $\lbrack \F(X) \rbrack$ is natural with respect to inclusions of $G$-spaces, which follows from the description using the projection.
\end{Rmk}

Using this we can now define the claimed identification as the map induced on colimits of the composite
\[\xymatrix@C=1.8cm{KK_*^G(C_0(X),\C) \ar[r]^-{\GJ} & KK_*(C_0(X)\rtimes G,\C) \ar[r]^-{-\circ \lbrack \F(X) \rbrack} & KK_*(C_0(X/G),\C).}\]

We want to conclude this construction by the following
\begin{Lemma}\label{diagram5}
Let $X$ be a proper, free and cocompact $G$-space. Then the inclusion \\$i: \C \to C(X/G)$ has the property
\[ i^*\lbrack \F(X)\rbrack = \lbrack\, p_X \rbrack \in KK(\C, C_0(X)\rtimes G).\]
\end{Lemma}
\begin{pf}
This follows immediately from the description of $\F(X)$ using the projection $p_X$.
\end{pf}

So we want to show that for each proper and cocompact $G$-space $X$ the composite
\begin{align}\label{composite}
\xymatrix@C=.7cm{KK^G_*(C_0(X),\C) \ar[r]^-{\GJ} & KK_*(C_0(X)\rtimes G, \C) \ar[r]^-{-\circ \lbrack \F(X) \rbrack} & KK_*(C(X/G),\C) \ar[r]^-{\MF} & KK_*(\C,C^*G)}
\end{align}
equals the analytical assembly map of section $2$. We recall that $\MF$ is defined by taking cup-cap product with the element 
\[ \lbrack \mathcal{L}_{X/G} \rbrack \in KK(\C,C(X/G)\otimes C^*G).\]
Now it is a standard fact from $KK$-theory, see for instance \cite[Proposition 18.9.1 (c)]{BB}, that the diagram
\[\xymatrix{KK_*(C_0(X)\rtimes G,\C) \ar[r]^-{\tau_{C^*G}} \ar[d]_{-\circ\lbrack \F(X) \rbrack} & KK_*((C_0(X)\rtimes G)\otimes C^*G,C^*G) \ar[d]^{-\circ \tau_{C^*G}\lbrack \F(X) \rbrack} \\
		KK_*(C(X/G),\C) \ar[r]_-{\tau_{C^*G}} & KK_*(C(X/G)\otimes C^*G,C^*G)}\]
commutes, which implies that the composite (\ref{composite}) is equal to
\[\xymatrix{KK^G_*(C_0(X),\C) \ar[r]^-{\GJ} & KK_*(C_0(X)\rtimes G, \C) \ar[r] & KK_*(\C,C^*G)}\]
where the last map is the cup-cap product with the class
\[ \tau_{C^*G}\lbrack \F(X)\rbrack \circ \lbrack \mathcal{L}_{X/G} \rbrack  \in KK(\C,(C_0(X)\rtimes G)
\otimes C^*G).\]
So our next goal is to factor the analytical assembly map
\[\xymatrix{KK^G_*(C_0(X),\C) \ar[r]^-{\GJ} & KK_*(C_0(X)\rtimes G,\C) \ar[r] & KK_*(\C,C^*G) }\]
in which the last map is cup-cap product with an element \[\beta_X \in KK(\C,(C_0(X)\rtimes G) \otimes C^*G).\]

The main to for this is the following Proposition. We want to thank Nigel Higson for pointing this out to us.

\begin{Prop}\label{factorization}
Let $A$ be a $G$-$C^*$-algebra. Then the following diagram commutes
\[\xymatrix@!C{KK_*^G(A,\C) \ar[r]^-{j^G} \ar[d]_{\GJ} & KK_*(A\rtimes G,C^*G) \\
KK_*(A\rtimes G,\C) \ar[r]_-{\tau_{C^*G}} & KK_*(A\rtimes G\otimes C^*G,C^*G) \ar[u]_{\Delta^*}.}\]
The similar statement for the reduced descent homomorphism is also true.
\end{Prop}
\begin{pf}
By classical results as for example in \cite{HR} we may assume that any element in $KK^G_*(A,\C)$ is represented by a triple $\lbrack \h,\pi,\F \rbrack$ where $\h$ is a separable Hilbert-space with unitary $G$-action $U: G \to \mathcal{B}(\h)$, $\pi:A \to \mathcal{B}(\h)$ is an equivariant representation and $\F \in \mathcal{B}(\h)$ is a selfadjoint operator satisfying the usual compactness conditions.

Let us first calculate what the lower composite does on such an element.
The Green-Julg map takes this class to the class $\lbrack \h, \pi\rtimes U, \F \rbrack$ where $\pi\rtimes U: A \rtimes G \to \mathcal{B}(\h)$ is induced by the covariant pair $(\pi,U)$. By definition we get that \[\tau_{C^* G}\lbrack \h,\pi\rtimes U,\F \rbrack = \lbrack \h\otimes C^*G, (\pi\rtimes U) \otimes \lambda, \F\otimes \id \rbrack.\]
where again $\lambda$ denotes the action by left multiplication of $C^*_rG$ on itself. Hence, we have 
\[ \Delta^*(\tau_{C^* G}(\mathrm{GJ}\lbrack \h,\pi, \F \rbrack )) =
\lbrack \h\otimes C^* G,((\pi\rtimes U)\otimes \lambda)\circ \Delta,\F\otimes \id \rbrack.\]
We need to compare this to $j^G\lbrack \h,\pi,\F \rbrack = \lbrack \h\rtimes G, \tilde{\pi},\tilde{\F} \rbrack$ and for this we will show the following facts:

\begin{enumerate}
\item[(i)] The Hilbert-$C^* G$-modules $\h\rtimes G$ and $\h\otimes C^* G$ are isomorphic and 
\item[(ii)] under this isomorphism, the operator $\tilde{\F}$ translates to $\F \otimes \id$ and the representation $\tilde{\pi}$ corresponds to $((\pi\rtimes U) \otimes \lambda)\circ\Delta$.
\end{enumerate}

For $\mathrm{(i)}$ we will begin by showing that $\h \otimes C_c(G,\C)$ and $C_c(G,\h)$ are isomorphic as $C_c(G,\C)$ modules. Note that it is clear that they are isomorphic as $\C$-modules so it suffices to check whether the canonical map $\h \otimes C_c(G,\C) \to C_c(G,\h)$ sending $x\otimes \alpha$ to the function $\alpha_x(g) = \alpha(g)\cdot x$ is a $C_c(G,\C)$ module map, which is a tedious but simple calculation.

Next we show that the given inner products on $C_c(G,\h)$ and $\h \otimes C_c(G,\C)$ coincide. More precisely let $\alpha, \beta \in C_c(G,\C)$ and $x,y \in \h$. Then as before we can view $\alpha_x$ and $\beta_y$ as elements of $C_c(G,\h)$. One can compute that
\[\langle \alpha_x,\beta_y \rangle_{C_c(G,\h)} = (\alpha^*\beta)\cdot \langle x,y \rangle_{\h} = \langle x\otimes \alpha, y\otimes\beta\rangle_{\h\otimes C^* G} .\]
It follows that there is an induced isomorphism of the Hilbert-$C^* G$-modules $\h\rtimes G$ and $\h\otimes C^* G$ as claimed.
\\ \\
In order to show $\mathrm{(ii)}$ we want that under this isomorphism $\tilde{\F}$ corresponds to $\F\otimes \mathrm{id}$.
This just means that for any $\alpha \in C_c(G,\C)$ and $x\in \h$ we want that
$\tilde{\F}(\alpha_x) = \alpha_{\F(x)}$, which is true since for any $g\in G$ we have
\begin{align*}
\tilde{\F}(\alpha_x)(g) = \F(\alpha_x(g)) = \F(\alpha(g)x) = \alpha(g)\F(x) = \alpha_{\F(x)}(g)
\end{align*}
as desired.

So it remains to show that the two representations
\[\xymatrix@R=.2cm@C=2.5cm{ A\rtimes G \ar[r]^-{\bar{\pi}\rtimes \bar{U}} & \ad(\h\rtimes G) \text{ and } \\
 A\rtimes G \ar[r]^-{((\pi\rtimes U)\otimes \lambda) \circ \Delta} & \ad(\h\otimes C^*G) }\]
correspond to each other under the canonical isomorphism of $\mathrm{(i)}$.
So let us calculate both representations.
Let $a \in A$ and $h \in G$, let $x\in \h$ and $\alpha \in C_c(G,\C)$. We define elements $\delta_h^a \in C_c(G,A)$ by the formula
\[ \delta^a_h(g) = \begin{cases} a & \text{ if } h=g, \\ 0 & \text{ else.} \end{cases}\]
It is the definition of $\Delta$ that we have
\[\Delta(\delta_h^a) = \delta_h^a \otimes h.\]

Hence we can compute
\begin{align*}
\left(((\pi\rtimes U) \otimes \lambda) \circ \Delta\right)(\delta^a_h)(x\otimes \alpha)(g) & = \left((\pi\rtimes U)\otimes \lambda\right)(\delta^a_h \otimes h)(x\otimes \alpha)(g) \\ & = \left((\pi\rtimes U)(\delta^a_h)(x)\otimes \lambda_h(\alpha)\right)(g) \\ & = \alpha(h^{-1}g)\cdot(\pi\rtimes U)(\delta^a_h)(x)  \\ & =
\alpha(h^{-1}g)\cdot\sum_{h^{\prime} \in G} \pi((\delta^a_h)(h^{\prime}))(U_{h^{\prime}}(x)) \\ & = \alpha(h^{-1}g) \cdot \pi(a)(U_h(x)).
\end{align*}
On the other hand by definition of the reduced descent homomorphism we can compute
\begin{align*}
(\tilde{\pi}\rtimes\tilde{U})(\delta^a_h)(\alpha_x)(g) & = \sum_{h^{\prime} \in G} \pi((\delta^a_h)(h^{\prime}))(U_{h^{\prime}}(\alpha_x(h^{\prime -1}g))) \\ & = \pi(a)(U_h(\alpha(h^{-1}g)x)) = \alpha(h^{-1}g)\cdot\pi(a)(U_h(x)).
\end{align*}
This completes the proof of the proposition.
\end{pf}
\begin{Rmk}
In \cite[section 2.2]{JR} Rosenberg claims that the (unreduced) descent homomorphism may be factored in a different way, but this factorization is not true. For example his factorizations says that for finite groups the diagram
\[\xymatrix{KK^G(\C,\C) \ar[r] \ar[d]_\GJ \ar[dr]|{j^G} & KK(\C,\C G) \ar[d]^{\varepsilon^*} \\ KK(\C G,\C) \ar[r]_{i_*} & KK(\C G,\C G)}\]
commutes, where the top horizontal map is the Green-Julg isomorphism for compact groups. Using that $j^G\lbrack \id_\C \rbrack = \lbrack \id_{\C G} \rbrack$ and that $\varepsilon\circ i = \id_\C$ this in turn implies the existence of elements $x \in KK(\C G,\C)$ and $y \in KK(\C,\C G)$ such that the identity map of $K_0(\C G)$ factors over $K_0(\C)$, hence $G$ would be the trivial group.

\end{Rmk}

The last proposition provides a factorization of the analytical assembly map as follows
\[\xymatrix{ & & KK_*(\C,C^*G) \\ KK_*(C_0(X)\rtimes G,\C) \ar@/^1pc/[urr] \ar[r]^-{\tau_{C^*G}} & KK_*((C_0(X)\rtimes G) \otimes C^*G,C^*G) \ar[r]^-{\Delta^*} & KK_*(C_0(X)\rtimes G, C^*G) \ar[u]_{-\circ \lbrack\, p_x \rbrack} \\
KK^G_*(C_0(X),\C) \ar@/_1pc/[urr]_-{j^G} \ar[u]^\GJ_\cong & &}\]
where now by definition of the cup-cap product the upper composite is just cup-cap product with the element
\[ \beta_X = \lbrack \Delta(p_X) \rbrack \in KK(\C,(C_0(X)\rtimes G) \otimes C^*G).\]
It is hence natural to ask whether we have an equality
\[ \lbrack \Delta(p_X) \rbrack = \tau_{C^*G}\lbrack \F(X) \rbrack \circ \lbrack \mathcal{L}_{X/G} \rbrack \in KK(\C,(C_0(X)\rtimes G)\otimes C^*G)\]
which would directly imply the main theorem. We recall that we have (using Lemma \ref{diagram5})
\[ \lbrack \Delta(p_X) \rbrack = \lbrack \Delta \rbrack \circ \lbrack\, p_X \rbrack = \lbrack \Delta \rbrack \circ i^*\lbrack \F(X) \rbrack.\]

The rest of this paper is devoted to a proof of this equality
\[ \lbrack \Delta \rbrack \circ i^*\lbrack \F(X) \rbrack = \tau_{C^*G}\lbrack \F(X) \rbrack \circ \lbrack \mathcal{L}_{X/G} \rbrack \in KK(\C,(C_0(X)\rtimes G)\otimes C^*G).\]

To do this we will apply methods from fixed point algebras as introduced by Kasparov in \cite[section 3]{KA} and more general versions as used in \cite{BE}. Once in the picture of fixed point algebras we can use results by Buss and Echterhoff in \cite{BE} to prove the equality. One of the crucial points in the proof is to relate the dual coaction $\Delta$ to the Hilbert-module that occurs in the Kasparov product $\tau_{C^*G}\lbrack \F(X) \rbrack \circ \lbrack \mathcal{L}_{X/G} \rbrack$.

\begin{Def}\label{B,X-modules}
Let $B$ be a $G$-$C^*$-algebra and $\e$ a $G$-Hilbert-$B$-module. If $\e$ is equipped with a $G$-equivariant morphism $C_0(X) \to\ad(\e)$ we call this datum a $(B,X\rtimes G)$\emph{-Hilbert-module}.
\end{Def}

\begin{Def}\label{G-objects}
We consider the following two $G$-$C^*$-algebras given by
\[ A = (C_0(X)\otimes C^*G,\tau\otimes\adl) \text{ and } \]
\[ B = (C_0(X)\otimes C^*G,\tau\otimes\id) \]
where $\adl$ denotes the conjugation action of $G$ on $C^*G$, and $\tau$ is the induced action of $G$ on $C_0(X)$. The object 
\[ \e = (C_0(X)\otimes C^*G,\tau\otimes\lambda) \]
naturally becomes an equivariant imprimitivity $A$-$B$-bimodule, where $\lambda$ denotes the left regular representation of $G$ on $C^*G$.
\end{Def}
We note that both $\e$ and $B$ are examples of $(B,X\rtimes G)$-Hilbert-modules as in Definition \ref{B,X-modules}, where the action $C_0(X) \to \ad(\e)$ and $C_0(X)\to \ad(B)$ is given by left multiplication on $C_0(X)$.

For the next definition see also the remark after \cite[Lemma 2.1]{EE} and the references listed there.
\begin{Def}\label{fixed point}
The \emph{generalized fixed point algebras} $A^G$ and $B^G$ are defined by:
\[ A^G = C_0(X\times_{G,\adl} C^*G) \text{ and }\]
\[ B^G = C_0(X\times_{G,\id} C^*G) \cong C_0(X/G)\otimes C^*G.\]
The \emph{generalized fixed point module} $\e^G$ is defined similarly by
\[ \e^G = C_0(X\times_{G,\lambda} C^*G).\]
\end{Def}
\begin{Rmk}
These algebras and this module may of course be interpreted as the algebras and the module of sections of the obvious bundles over $X/G$. In this notation we have that $\e^G = \Gamma_0(\mathcal{L}_{X/G})$.
\end{Rmk}

\begin{Lemma}
The algebras $A^G$ and $B^G$ are Morita equivalent.
\end{Lemma}
\begin{pf}
The generalized fixed point module $\e^G$ is an imprimitivity bimodule.
\end{pf}
\begin{Rmk}
We just want to emphasize again that this implies that the element
\[ \lbrack\, \e^G \rbrack \in KK(A^G,B^G) \]
is a $KK$-equivalence.
\end{Rmk}

There exists an inclusion $j:C_0(X/G) \to A^G$ using the fact that the conjugation action of $G$ on $\C \subset C^*G$ is trivial. Using that $B^G = C_0(X/G)\otimes C^*G$ we have the following

\begin{Lemma}\label{diagram4}
If $X$ is in addition cocompact we have
\[ \lbrack \e^G,\pi_\C,0 \rbrack = \lbrack \mathcal{L}_{X/G} \rbrack \in KK(\C,C(X/G)\otimes C^*G)\]
where $\pi_\C$ is the unique unital representation.
\end{Lemma}
\begin{pf}
This follows directly from the definitions. 
\end{pf}

Following \cite[section 3]{BE}, we need to extend the construction of $\F(X)$ to a more general situation.
\begin{Def}\label{F general}
Given any $(B,X\rtimes G)$-Hilbert-module $\e$ we define
\[ \F_c(\e) = C_c(X)\cdot \e \]
which can be viewed as a right $C_c(G,B)$-module and as left $C_0(X/G)$-module. Moreover $\F_c(\e)$ has a $C_c(G,B)$-valued inner product, with respect to which it completes to a Hilbert-$B\rtimes G$-module $\F(\e)$.
\end{Def}

\begin{ex}
Consider $B = C_0(X) = \e$ as a $G$-Hilbert-$C_0(X)$-module over itself. The action by multiplication operators is $G$-equivariant and so we get a Hilbert-$C_0(X)\rtimes G$-module $\F(C_0(X))$ and it can be checked that 
\[ \F(C_0(X)) \cong \F(X).\]
\end{ex}

We need the following technical observations.

\begin{Lemma}\label{identification}
In the notation of Definition \ref{G-objects} we have
$\e \rtimes G = B\rtimes G = (C_0(X)\rtimes G)\otimes C^*G$.
Moreover the multiplication action of $C_0(X)$ on $\e$ is equivariant and the induced action
\[\xymatrix{C_0(X) \rtimes G \ar[r] & \ad(\e\rtimes G) = \ad(C_0(X)\rtimes G \otimes C^*G) }\]
may be identified with the multiplication action after applying the dual coaction. Precisely the diagram
\[\xymatrix{C_0(X)\rtimes G \ar[r] \ar[d]_{\Delta} & \ad(C_0(X)\rtimes G \otimes C^*G) \\
C_0(X)\rtimes G \otimes C^*G \ar@/_1.2pc/[ur]_-{M} &}\]
commutes, where $\Delta$ is the dual coaction and $M$ is the left-multiplication action.
Furthermore $\F(B) \cong \F(X)\otimes C^*G$ as right-$B\rtimes G$-modules.
\end{Lemma}
\begin{pf}
It is clear that $\e\rtimes G = B\rtimes G$ because the $G$-action on $\e$ is not used when constructing $\e\rtimes G$ (only the action on $B$ is relevant for this). The $G$-action is used when constructing the action map
$C_0(X)\rtimes G \to \ad(\e\rtimes G)$. We recall from \cite[section 2.2, formula 2.5]{BE} that $C_c(G,C_0(X))$ acts on $C_c(G,\e)$ by the formula
\[ (f.\alpha)(g,x) = \sum_{h \in G} f(h,x)\cdot h \cdot \alpha(h^{-1}g,h^{-1}x) \]
for $f \in C_c(G,C_0(X))$ and $\alpha \in C_c(G,C_0(X)\otimes C^*G)$. Note that the extra $h$ in the product comes precisely from the $G$-action $\tau\otimes\lambda$ on $\e$.
Now by definition of the convolution product on $B\rtimes G$ we have that
\[ (\Delta(f)\ast \alpha)(g,x) = \sum_{h\in G} \Delta(f)(h,x)\cdot\alpha(h^{-1}g,h^{-1}x).\]
Here, no extra $h$-factor comes up in the product with $\alpha$ as the $G$-action on $B$ is given by $\tau\otimes\id$, i.e. is trivial on the $C^*G$-term.

Now we consider the special functions $f = \delta_h^\varphi \in C_c(G,C_0(X))$ for $\varphi \in C_0(X)$. We recall that these are given by
\[ \delta^\varphi_h(g) = \begin{cases} \varphi & \text{ if } h=g, \\ 0 & \text{ else,} \end{cases}\]
and as in Proposition \ref{factorization} we have
\[\Delta(\delta_h^\varphi) = \delta_h^\varphi \otimes h.\]
Using this we can easily see that
\[ (\delta^\varphi_h . \alpha)(g,x) = \varphi(x)\cdot h\cdot \alpha(h^{-1}g,h^{-1}x) = (\Delta(\delta^\varphi_h)\ast\alpha)(g,x) \]
which shows the commutativity of the diagram.

The statement about $\F(B)$ follows directly from the fact that the $G$-action is trivial on the $C^*G$-tensor factor and the example previous to this lemma.
\end{pf}

\begin{Prop}\label{BE 3.6}
In the situation of Definition \ref{F general} there is an isomorphism of Hilbert-$B\rtimes G$-modules
\[\xymatrix{\Psi: \F(X) \otimes_{C_0(X)\rtimes G} (\e\rtimes G) \ar[r]^-{\cong} & \F(\e)}\]
where $\e\rtimes G$ is as in Lemma \ref{descent}. 
\end{Prop}
\begin{pf}
This is a special case of \cite[Proposition 3.6]{BE}. 
\end{pf}

It turns out (using the description in \cite[before Prop. 3.20]{BE}) that the equivariant fixed point module $\e^G$ as defined in Definition \ref{fixed point} coincides with the more general construction as in \cite[Lemma 4.1]{BE}. In particular we have the following characterization of $\e^G$ in terms of the construction of Definition \ref{F general}.

\begin{Prop}\label{BE 4.6}
There is an isomorphism of Hilbert-$B^G$-modules
\[\xymatrix{\e^G \ar[r]^-{\cong} & \F(\e) \otimes_{B\rtimes G} \F(B)^*}.\]

\end{Prop}
\begin{pf}
This is proven in \cite[Corollary 4.6]{BE} using \cite[Proposition 4.5]{BE}.
\end{pf}

\noindent

We have now collected all results needed to prove the remaining equality. So let us start by computing the Kasparov product 
\[\tau_{C^*G}\lbrack \F(X) \rbrack \circ \lbrack \mathcal{L}_{X/G} \rbrack = \tau_{C^*G}\lbrack \F(X) \rbrack \circ \lbrack \,\e^G,\pi_\C,0 \rbrack.\]
First, we claim that
\[ \tau_{C^*G}\lbrack \F(X) \rbrack = \lbrack \F(B) \rbrack \in KK(B^G,B\rtimes G).\]
Indeed, by Lemma \ref{identification} we have that $\F(B) \cong \F(X)\otimes C^*G$ and $B^G \cong C_0(X/G)\otimes C^* G$ and the left $B^G$-module action on $\F(B)$ corresponds precisely to \[\pi\otimes \l : C_0(X/G) \otimes C^*G \to \ad(\F(X)\otimes C^*G)\]
as in the definition of the exterior product $\tau_{C^*G}$.
Hence we can compute the Kasparov product to be
\[ \lbrack \F(B) \rbrack \circ \lbrack \e^G,\pi_\C,0 \rbrack = \lbrack \e^G\otimes_{B^G} \F(B), \pi_\C, 0 \rbrack \in KK(\C,B\rtimes G).\]
The fact that this is a Kasparov product follows from the construction of it (the connection one needs to construct may be chosen to be zero in this case, as all operators involved are the zero operators).

Altogether this means that
\[ \tau_{C^*G}\lbrack \F(X) \rbrack \circ \lbrack \mathcal{L}_{X/G} \rbrack = \lbrack \e^G\otimes_{B^G} \F(B), \pi_\C, 0 \rbrack \in KK(\C,B\rtimes G).\]
We can now compute the Hilbert-$B\rtimes G$-module $\e^G\otimes_{B^G} \F(B)$ as follows

\begin{align*}
\e^G\otimes_{B^G} \F(B) & \cong \F(\e)\otimes_{B\rtimes G} \F(B)^* \otimes_{B^G} \F(B) & \text{ by Prop }\ref{BE 4.6}\\
& \cong  \F(\e) & \\
& \cong \F(X)\otimes_{C_0(X)\rtimes G} \left(\e\rtimes G \right) &\text{ by Prop } \,\;\ref{BE 3.6} \\
& \cong \F(X) \otimes_{\Delta} B\rtimes G & \text{ by Prop }\,\;\ref{identification}
\end{align*}
Now on the other hand we want to compute the element
\[ \lbrack \Delta \rbrack \circ i^*\lbrack \F(X) \rbrack  \in KK(\C,B\rtimes G)\]
but by sheer definition we get that
\[ \lbrack \Delta \rbrack \circ i^*\lbrack \F(X) \rbrack = \lbrack \F(X)\otimes_{\Delta} B\rtimes G, \pi_\C, 0 \rbrack.\]
Hence we have
\begin{align*}
\tau_{C^*G}\lbrack \F(X) \rbrack \circ \lbrack \mathcal{L}_{X/G} \rbrack & = \lbrack \e^G \otimes_{B^G} \F(B),\pi_\C,0 \rbrack \\ & = \lbrack \F(X)\otimes_\Delta B\rtimes G,\pi_\C,0 \rbrack \\ & = \lbrack \Delta \rbrack \circ i^*\lbrack \F(X) \rbrack.
\end{align*}

\end{section}

\begin{Rmk}
At last we want to point out that our main theorem implies that the usual Baum-Connes assembly map (which is the reduced analytical assembly map) is also equal to a Mishchenko-Fomenko index if we replace $C^*G$ by $C^*_r G$ throughout the whole paper.
\end{Rmk}

\end{document}